\documentclass[12pt]{article}
\usepackage{amsmath}
\usepackage{amssymb}
\usepackage{amsthm}
\usepackage{amscd}
\usepackage{amsfonts}
\usepackage{enumerate}
\usepackage[hidelinks]{hyperref}
\newtheorem{teo}{Theorem}

\begin{document}
\title{Some Open Problems Concerning Orthogonal Polynomials on Fractals and Related Questions}
\author{G\"{o}kalp Alpan and Alexander Goncharov}
\date{}
\maketitle

\section{Background and notation}
\subsection{Chebyshev and orthogonal polynomials}

Let $K\subset\mathbb{C}$ be a compact set containing infinitely many points.
We use $\|\cdot\|_{L^\infty(K)}$ to denote the sup-norm on $K$,
$\mathcal{M}_n$ is the set of all monic polynomials of degree $n$. The polynomial $T_{n,K}$
that minimizes $\|Q_n\|_{L^\infty(K)}$ for $Q_n \in \mathcal{M}_n$ is called the $n$-th \emph{Chebyshev polynomial} on $K$.

Let the  logarithmic capacity $\mathrm{Cap}(K)$ be positive. Then we define the $n$-th Widom factor for $K$ by $$W_{n}(K):= \|T_{n,K}\|_{L^\infty(K)}/\mathrm{Cap}(K)^n.$$

In what follows we consider unit Borel measures $\mu$ with non-polar compact support $\mathrm{supp}(\mu)$ in $\mathbb{C}$. The $n$-th monic \emph{orthogonal polynomial} $P_n(z;\mu) = z^n + \ldots$
associated with $\mu$ has the property
$$\displaystyle \|P_n(\cdot;\mu)\|_ {L^2(\mu)}^2= \inf_{Q_n \in \mathcal{M}_n}\int |Q_n(z)|^2\,d\mu(z),$$
where $\| \cdot \|_{L^2(\mu)}$ is the norm in $L^2(\mu)$. Then the $n$-th \emph{Widom-Hilbert factor} for $\mu$ is $$W_n^2(\mu):=\|P_n(\cdot;\mu)\|_{L^2(\mu)}/(\mathrm{Cap}(\mathrm{supp}(\mu)))^n.$$

If $\mathrm{supp(\mu)}\subset\mathbb{R}$ then a three-term recurrence relation
$$ x P_n(x;\mu)= P_{n+1}(x;\mu)+b_{n+1} P_n(x;\mu)+a_n^2 P_{n-1}(x;\mu)$$
is valid for $n\in \mathbb{N}_0:=\mathbb{N}\cup\{0\}.$ The initial conditions $P_{-1}(x;\mu)\equiv0$ and $P_{0}(x;\mu)\equiv1$ generate two bounded sequences $(a_n)_{n=1}^\infty, (b_n)_{n=1}^\infty$ of \emph{recurrence coefficients} associated with $\mu$. Here, $a_n>0,\,b_n\in \mathbb{R}$ for  $n\in\mathbb{N}$ and $$ \|P_n(\cdot;\mu)\|_ {L^2(\mu)}=a_1\cdots a_n.$$

A bounded  two sided $\mathbb{C}$-valued sequence  $(d_n)_{n=-\infty}^\infty$ is called \emph{almost periodic} if the set $\{(d_{n+k})_{n=-\infty}^\infty:\,\, k\in\mathbb{Z}\}$ is precompact in $l^\infty(\mathbb{Z})$. A one sided sequence $(c_n)_{n=1}^\infty$ is called almost periodic if it is the restriction of a two sided almost periodic sequence to $\mathbb{N}$. A sequence $(e_n)_{n=1}^\infty$ is called \emph{asymptotically almost periodic} if there is an almost periodic sequence $(e^\prime_n)_{n=1}^\infty$ such that $|e_n-e^\prime_{n}|\rightarrow 0$ as $n\rightarrow 0$.

A class of Parreau-Widom sets plays a special role in the recent theory of orthogonal and Chebyshev polynomials.
Let $K$ be a non-polar compact set and $g_{\mathbb{C}\setminus K}$ denote the Green function for $\overline{\mathbb{C}}\setminus K$ with a pole at infinity. Suppose $K$ is regular with respect to the Dirichlet problem, so the set $\mathcal{C}$ of critical points of $g_{\mathbb{C}\setminus K}$ is at most countable.
Then $K$ is said to be a \emph{Parreau-Widom} set if $\sum_{c\in \mathcal{C} } g_{\mathbb{C}\setminus K}(c)<\infty$.
Parreau-Widom sets on $\mathbb{R}$ have positive Lebesgue measure. For different aspects of such sets, see \cite{christiansen,hasumi,yuditskii}.\\

A class of regular in the sense of Stahl-Totik measures can be defined by the following condition
$$\lim_{n\rightarrow\infty}W_n(\mu)^{1/n}=1.$$\\

For a measure $\mu$ supported on $\mathbb{R}$ we use the  Lebesgue decomposition of $\mu$ with respect to the
Lebesgue measure:
$$d\mu(x)=f(x)dx+d\mu_{\mathrm{s}}(x).$$

Following \cite{christi}, let us define the Szeg\H{o} class $\mathrm{Sz}(K)$ of measures on a given Parreau-Widom set $K\subset \mathbb{R}.$ Let $\mu_K$ be the equilibrium measure on $K$.
By $\mathrm{ess\,supp}(\cdot)$ we denote the essential support of the measure, that is the set of accumulation points of the support. We have $\mathrm{Cap(supp}(\mu))= \mathrm{Cap(ess\,supp}(\mu)),$ see Section 1 of \cite{Sim3}. A measure $\mu$ is in the \emph{Szeg\H{o} class} of $K$ if

\begin{enumerate}[(i)]
	\item $\mathrm{ess\,supp}(\mu)=K.$
	\item $\int_K \log{f(x)}\,d\mu_K(x)>-\infty.$ (Szeg\H{o} condition)
	\item the isolated points $\{x_n\}$ of $\mathrm{supp}(\mu)$
	satisfy $\sum_n g_{\mathbb{C}\setminus K}(x_n)<\infty.$
\end{enumerate}
By Theorem 2 in \cite{christi} and its proof, (ii) can be replaced by one of the following conditions:
\begin{enumerate}[(ii$^\prime$)]
	\item $\limsup_{n\rightarrow\infty}W_n^2(\mu)>0.$ (Widom condition)
\end{enumerate}
\begin{enumerate}[(ii$^{\prime\prime}$)]
	\item $\liminf_{n\rightarrow\infty}W_n^2(\mu)>0.$ (Widom condition 2)
\end{enumerate}

 One can show that any $\mu\in\mathrm{Sz}(K)$ is regular in the sense of Stahl-Totik.

\subsection{Generalized Julia sets and $K(\gamma)$}
Let $(f_n)_{n=1}^\infty$ be a sequence of rational functions with $\deg{f_n}\geq 2$ in $\overline{\mathbb{C}}$ and $F_n:= f_n \circ f_{n-1}\circ\ldots \circ f_{1}$. The domain of normality for $(F_n)_{n=1}^\infty$ in the sense of Montel is called the \emph{Fatou set} for $(f_n)$. The complement of the Fatou set in $\overline{\mathbb{C}}$ is called the \emph{Julia set} for $(f_n)$. We denote them by $F_{(f_n)}$ and $J_{(f_n)}$ respectively. These sets were considered first in \cite{Fornaess}. In particular, if $f_n=f$ for some fixed rational function $f$ for all $n$ then  $F_{(f)}$ and $J_{(f)}$ are used instead. To distinguish this last case, the word \emph{autonomous} is used in the literature.

Suppose $f_n(z)=\sum_{j=0}^{d_n}a_{n,j}\cdot z^j$ where $d_n \geq 2$ and $a_{n,d_n}\neq 0$ for all $n\in\mathbb{N}$.
Following \cite{Bruck}, we say that $(f_n)$ is a \emph{regular polynomial sequence} (write $(f_n)\in\mathcal{R}$) if positive constants
$A_1, A_2, A_3$ exist such that for all $n\in\mathbb{N}$ we have the following three conditions:\\
$|a_{n,d_n}|\geq A_1$\\
$|a_{n,j}|\leq A_2 |a_{n,d_n}|$ for $j=0,1,\ldots, d_n-1$\\
$\log{|a_{n,d_n}|}\leq A_3\cdot d_n$\\
For such polynomial sequences, by \cite{Bruck}, $J_{(f_n)}$ is a regular compact set in  $\mathbb{C}$.
In addition, $\mathrm{Cap}(J_{(f_n)})>0$ and $J_{(f_n)}$ is the boundary of $$\mathcal{A}_{(f_n)}(\infty):=\{z\in\overline{\mathbb{C}}: F_n(z) \mbox{ goes locally uniformly to } \infty \}.$$

The following construction is from \cite{gonc}. Let $\gamma:=(\gamma_k)_{k=1}^\infty$ be a sequence provided that $0<\gamma_k<1/4$ holds for all $k\in\mathbb{N}$ and $\gamma_0:=1$. Let $f_1(z)=2z(z-1)/\gamma_1+1$ and $f_n(z)=\frac{1}{2\gamma_n}(z^2-1)+1$ for $n>1$. Then $K(\gamma):=\cap_{s=1}^{\infty} F_s^{-1}([-1,1])$ is a Cantor set on $\mathbb{R}$. Furthermore, $F_s^{-1}([-1,1])\subset F_t^{-1}([-1,1])\subset [0,1]$ whenever $s>t$.

Also we use an expanded version of this set. For a sequence $\gamma$ as above, let $f_n(z)=\frac{1}{2\gamma_n}(z^2-1)+1$ for $n\in\mathbb{N}$. Then $K_1(\gamma):=\cap_{s=1}^{\infty} F_s^{-1}([-1,1])\subset [-1,1]$ and $F_s^{-1}([-1,1])\subset F_t^{-1}([-1,1])\subset [-1,1]$ provided that $s>t$. It is a Cantor set. If there is a $c$ with $0<c<\gamma_k$ for all $k$ then $(f_n)\in\mathcal{R}$ and $J_{(f_n)}=K_1(\gamma)$, see \cite{g5}. If $\gamma_1=\ldots=\gamma_k$ for all $k\in\mathbb{N}$ then $K_1(\gamma)$ is an autonomous polynomial Julia set.

\subsection{Hausdorff measure}
A function $h:{\mathbb{R}}_+ \to {\mathbb{R}}_+$ is called a \emph{dimension function} if it is increasing, continuous
and $h(0)=0.$  Given a set $E\subset \mathbb{C},$ its \emph{$h$-Hausdorff measure} is defined as
\begin{equation*}
\Lambda_h(E)=\lim_{\delta \rightarrow 0} \inf \left\{ \sum h(r_j) : E \subset \bigcup B(z_j, r_j)\,\,\,
\mbox {with}\,\,\,   r_j \leq \delta \right\},
\end{equation*}
where $B(z, r)$ is the open ball of radius $r$ centered at $z$. For a dimension function $h$,  a set $K\subset \mathbb{C}$ is an $h$-set if  $0<\Lambda_h(K) < \infty$. To denote the Hausdorff measure for $h(t)=t^\alpha$, $\Lambda_\alpha$ is used. \emph{Hausdorff dimension} of $K$ is defined as $\inf\{\alpha\geq 0: \Lambda_\alpha(K)=0\}$.

\section{Smoothness of Green functions and Markov Factors}

The next set of problems is concerned with the smoothness properties of the Green function
$g_{\, {\Bbb C} \setminus K}$ near compact set $K$ and related questions
We suppose that $K$ is regular with respect to the Dirichlet problem, so the function
$g_{\, {\Bbb C} \setminus K}$ is continuous throughout $\Bbb C.$ The next problem was posed in \cite{gonc}.\\

{\bf Problem 1.} Given modulus of continuity $\omega$, find a compact set $K$ such that the modulus of
continuity $\omega(g_{\, {\Bbb C}\setminus K}, \cdot)$ is similar to $\omega$.\\

Here, one can consider similarity as coincidence of moduli on some null sequence or in the sense of
weak equivalence: $\exists C_1, C_2$ such that
$$  C_1\,\omega(\delta) \leq \omega(g_{\, {\Bbb C}\setminus K}, \delta)\leq C_2\,\omega(\delta) $$
for sufficiently small positive $\delta$.

We guess that a set $K(\gamma)$ from \cite{gonc} is a candidate for the desired $K$ provided a suitable choice of
the parameters. We recall that, for many moduli of continuity, the corresponding Green's functions were
given in \cite{gonc}, whereas the characterization of optimal  smoothness for $g_{\,{\Bbb C} \setminus K(\gamma)}$
is presented in [\cite{g5}, Th.6.3]. \\

A stronger version of the problem is about pointwise estimation of the Green function:\\

{\bf Problem 2.} Given modulus of continuity $\omega$, find a compact set $K$ such that
$$ C_1 \, \omega(\delta) \leq g_{\,{\Bbb C}\setminus K}(z) \leq C_2 \, \omega(\delta)$$
for $\delta=dist(z,K)\leq \delta_0,$ where $C_1, C_1$ and $\delta_0$ do not depend on $z$. \\

In the most important case we get a problem about ``two-sided H{\"o}lder" Green function, which was posed by
A. Volberg on his seminar (quoted with permission):\\

{\bf Problem 3.} Find a compact set $K$ on the line such that for some $\alpha>0$ and constants $C_1, C_2$, if $\delta=dist(z,K)$ is small enough then
 \begin{equation}
C_1 \, \delta^{\alpha} \leq g_{\,{\Bbb C}\setminus K}(z) \leq C_2 \, \delta^{\alpha}.
 \end{equation}

Clearly, a closed analytic curve gives a solution for sets on the plane.

If $K \subset \Bbb R$ satisfies (1), then $K$ is of Cantor-type. Indeed, if
interior of $K$ (with respect to $\Bbb R$) is not empty, let $(a,b) \subset K,$ then
$g_{\,{\Bbb C}\setminus K}$ has $Lip\,1$ behavior near the point $(a+b)/2$. On the other hand,
near endpoints of $K$ the function $g_{\,{\Bbb C}\setminus K}$ cannot be better than $Lip\,1/2.$\\

By the Bernstein-Walsh inequality, smoothness properties of Green functions are closely related with
a character of maximal growth of polynomials outside the corresponding compact sets, which, in turn,
allows to evaluate Markov's factors for the sets. Recall that, for a fixed $n\in {\Bbb N}$ and (infinite)
compact set $K$, the $n-$th Markov factor $M_n(K)$ is the norm of operator of differentiation in the
space of holomorphic polynomials  ${\cal P}_n$ with the uniform norm on $K$. In particular, the H{\"o}lder
smoothness (the right inequality in (1)) implies Markov's property of the set $K$  (a polynomial growth rate
of $M_n(K)$). The problem about inverse implication (see e.g [20]) has attracted attention of many
researches.

By W. Ple\'sniak [20], any Markov set $K\subset {\Bbb R}^d$ has the extension property ($EP$), which means
that there exists  a continuous linear extension operator from the space of Whitney functions ${\mathcal E}(K)$
to  the space of infinitely differentiable functions on ${\Bbb R}^d.$ We guess that there is some extremal growth rate of $M_n$ which implies the lack of $EP$. Recently it was shown in \cite{gonczel}
that there is no complete characterization of $EP$ in terms of growth rate of Markov's factors. Namely,
two sets were presented, $K_1$ with $EP$ and $K_2$ without it, such that $M_n(K_1)$ grows essentially
faster than $M_n(K_2)$ as $n\to \infty$. Thus there exists non-empty zone of uncertainty
where the growth rate of $M_n(K)$ is not related with $EP$ of the set $K$.\\

{\bf Problem 4.} Characterize the growth rates of Markov's factors that define the boundaries of the zone of uncertainty for the extension property.

\section{Orthogonal polynomials}
One of the most interesting problems concerning orthogonal polynomials on Cantor sets on $\mathbb{R}$ is the character of periodicity of recurrence coefficients. It was conjectured in p. 123 of \cite{belis} that if $f$ is a non-linear polynomial such that $J(f)$ is a totally disconnected subset of $\mathbb{R}$ then the recurrence coefficients for $\mu_{J(f)}$ are almost periodic. This is still an open problem. In \cite{g8}, the authors conjectured that the recurrence coefficients for $\mu_{K(\gamma)}$ are asymptotically almost periodic for any $\gamma$. It may be hoped that a more general and slightly weaker version of Bellissard's conjecture can be valid.\\

{\bf Problem 5.} Let $(f_n)$ be a regular polynomial sequence such that $J_{(f_n)}$ is a Cantor-type subset of the real line. Prove that the recurrence coefficients for $\mu_{J(f_n)}$ are asymptotically almost periodic.\\

For a measure $\mu$ which is supported on $\mathbb{R}$, let $Z_n(\mu):=\{x:\,P_n(x;\mu)=0\}$. We define $U_n(\mu)$ by
\begin{equation*}
U_n(\mu):= \inf_{\substack{x,x^\prime\in Z_n(\mu)\\ x\neq x\prime}}|x-x^\prime|.
\end{equation*}

In \cite{kruger} Kr\"{u}ger and Simon gave a lower bound for $U_n(\mu)$ where $\mu$ is the Cantor-Lebesgue measure of the (translated and scaled) Cantor ternary set. In \cite{jons}, it was shown that Markov's inequality and spacing of the zeros of orthogonal polynomials are somewhat related.

Let $\gamma=(\gamma_k)_{k=1}^\infty$ and $n\in\mathbb{N}$ with $n>1$ be given and define $\delta_k=\gamma_0 \cdots \gamma_k$ for all $k\in\mathbb{N}_0$. Let $s$ be the integer satisfying $2^{s-1}\leq n < 2^s$. By \cite{g2},
\begin{equation*}\label{eq1}
\delta_{s+2}\leq U_n(\mu_{K(\gamma)})\leq \frac{\pi^2}{4}\cdot \delta_{s-2}
\end{equation*}
holds. In particular, if there is a number $c$ such that $0<c<\gamma_k<1/4$ holds for all $k\in\mathbb{N}$ then, by \cite{g2}, we have
\begin{equation}\label{eq2}
c^2\cdot \delta_{s}\leq U_n(\mu_{K(\gamma)})\leq \frac{\pi^2}{4c^2}\cdot \delta_{s}.
\end{equation}
By [13], at least for small sets $K(\gamma),$ we have $M_{2^s}(K(\gamma))\sim 2/\delta_s,$ where the symbol $\sim$
means the strong equivalence.\\

{\bf Problem 6.} Let $K$ be a non-polar compact subset of $\mathbb{R}$. Is there a general relation between the zero spacing of orthogonal polynomials for $\mu_{K}$ and smoothness of $g_{\mathbb{C}\setminus K}$? Is there a relation between the zero spacing of $\mu_K$ and the Markov factors?\\

As mentioned in section 1, the Szeg\H{o} condition and the Widom condition are equivalent for Parreau-Widom sets. Let $K$ be a Parreau-Widom set. Let $\mu$ be a measure such that $\mathrm{ess\,supp}(\mu)=K$ and the isolated points $\{x_n\}$ of $\mathrm{supp}(\mu)$
satisfy $\sum_n g_{\mathbb{C}\setminus K}(x_n)<\infty$. Then, as it is discussed in Section 6 of \cite{g4}, the Szeg\H{o} condition is equivalent to the condition
\begin{equation}\label{cond1}
\int_K \log(d\mu/d\mu_K)\,d\mu_K(x)>-\infty.
\end{equation}
This condition is also equivalent to the Widom condition under these assumptions.

It was shown in \cite{g9} that $\inf_{n\in\mathbb{N}} W_n(\mu_K)\geq 1$ for non-polar compact $K \subset \mathbb{R}.$
Thus the Szeg\H{o} condition in the form \eqref{cond1} and the Widom condition are related on arbitrary non-polar sets.\\

{\bf Problem 7.} Let $K$ be a non-polar compact subset of $\mathbb{R}$ which is regular with respect to the Dirichlet problem. Let $\mu$ be a measure such that $\mathrm{ess\,supp}(\mu)=K$. Assume that the isolated points $\{x_n\}$ of $\mathrm{supp}(\mu)$ satisfy $\sum_n g_{\mathbb{C}\setminus K}(x_n)<\infty$. If the condition \eqref{cond1} is valid for $\mu$ is it necessarily true that the Widom condition or the Widom condition 2 holds? Conversely, 
does the Widom condition imply \eqref{cond1}?\\

It was proved in \cite{csz} that if $K$ is a Parreau-Widom set which is a subset of $\mathbb{R}$ then $(W_n(K))_{n=1}^\infty$ is bounded above. On the other hand, $(W_n(K))_{n=1}^\infty$ is unbounded for some
Cantor-type sets, see e.g. \cite{gonchat}.\\

{\bf Problem 8.} Is it possible to find a regular non-polar compact subset $K$ of $\mathbb{R}$ which is not Parreau-Widom but $(W_n(K))_{n=1}^\infty$ is bounded? If $K$ has zero Lebesgue measure then is it true that $(W_n(K))_{n=1}^\infty$ is unbounded? We can ask the same problems if we replace $(W_n(K))_{n=1}^\infty$ by $(W_n^2(\mu_K))_{n=1}^\infty$ above.\\

Let $T_N$ be a real polynomial of degree $N$ with $N\geq 2$ such that it has $N$ real and simple zeros $x_1<\dots <x_n$ and $N-1$ critical points $y_1<\dots<y_{n-1}$ with $|T_N(y_i)|\geq 1$ for each $i\in\{1,\ldots,N-1\}$. We call such a polynomial admissible. If $K=T_N^{-1}([-1,1])$ for an admissible polynomial $T_N$  then $K$ is called a $T$-set. The following result is well known, see e.g. \cite{totikacta}.

\begin{teo}\label{thm1}
	Let $K=\cup_{j=1}^n [\alpha_j,\beta_j]$ be a disjoint union of $n$ intervals such that $\alpha_1$ is the leftmost end point. Then $K$ is a $T$-set if and only if $\mu_K([\alpha_1, c])$ is in $\mathbb{Q}$ for all $c\in\mathbb{R}\setminus K$.
\end{teo}

For $K(\gamma)$, it is known that $\mu_{K(\gamma)}([0,c])\in\mathbb{Q}$ if $c\in\mathbb{R}\setminus K(\gamma)$, see Section 4 in \cite{g2}.\\

{\bf Problem 9.} Let $K$ be a regular non-polar compact subset of $\mathbb{R}$ and $\alpha$ be the leftmost end point of $K$. Let $\mu_K([\alpha, c])\in \mathbb{Q}$ for all $c\in \mathbb{R}\setminus K$. What can we say about $K$? Is it necessarily a polynomial generalized Julia set? Does this imply that there is a sequence of admissible polynomials $(f_n)_{n=1}^\infty$ such that $(F_n^{-1}{[-1,1]})_{n=1}^\infty$ is a decreasing sequence of sets such that $K=\cap_{n=1}^\infty F_n^{-1}{[-1,1]}$?

\section{Hausdorff measures}

It is valid for a wide class of Cantor sets that the equilibrium measure and the corresponding Hausdorff measure on this set are mutually singular, see e.g. \cite{mak}.

Let $\gamma=(\gamma_k)_{k=1}^\infty$ with $0<\gamma_k<1/32$ satisfy $\sum_{k=1}^\infty \gamma_k <\infty$. This implies that $K(\gamma)$ has Hausdorff dimension $0$. In \cite{g3}, the authors constructed a dimension function $h_\gamma$ that makes $K(\gamma)$ an $h$-set. Provided also that $K(\gamma)$ is not polar it was shown that there is a $C>0$ such that for any Borel set $B$, $C^{-1}\cdot \mu_{K(\gamma)}(B)<\Lambda_{h_\gamma}(B)<C\cdot \mu_{K(\gamma)}(B)$ and in particular the equilibrium measure and $\Lambda_{h_\gamma}$ restricted to $K(\gamma)$ are mutually absolutely continuous. In \cite{gonczel}, it was shown by the authors that indeed these two measures coincide. To the best of our knowledge, this is the first example of a subset of $\mathbb{R}$ such that the equilibrium measure is a Hausdorff measure restricted to the set.\\

{\bf Problem 10.} Let $K$ be a non-polar compact subset of $\mathbb{R}$ such that $\mu_K$ is equal to a Hausdorff measure restricted to $K$. Is it necessarily true that the Hausdorff dimension of $K$ is $0$?\\

Hausdorff dimension of a unit Borel measure $\mu$ supported on $\mathbb{C}$ is defined by $\mathrm{dim}(\mu):=\inf \{\mathrm{HD}(K): \mu(K)=1\}$ where $HD(\cdot)$ denotes Hausdorff dimension of the given set. For polynomial Julia sets which are totally disconnected there is a formula for $\mathrm{dim}(\mu_{J(f)})$, see e.g.p. 23 in \cite{mak} and p.176-177 in \cite{prz}.\\

{\bf Problem 11.} Is it possible to find simple formulas for $\mathrm{dim}\left(\mu_{J_{(f_n)}}\right)$ where $(f_n)$ is a regular polynomial sequence?

\thanks{The authors are partially supported by a grant from T\"{u}bitak: 115F199}

\end{document}